\documentclass[11pt]{article}
\usepackage[T1]{fontenc}
\usepackage{lmodern,amsmath,amssymb,amsthm,mathtools,hyperref}
\usepackage[margin=1in]{geometry}
\newtheorem{theorem}{Theorem}[section]
\newtheorem{proposition}[theorem]{Proposition}
\newtheorem{remark}[theorem]{Remark}
\newtheorem{example}[theorem]{Example}
\newcommand{\HS}{\mathrm{HS}}
\newcommand{\Lop}{\mathfrak{L}}

\newcommand{\Pp}{\mathrm{P}}
\newcommand{\dd}{\,\mathrm d}
\newcommand{\norm}[1]{\left\|#1\right\|}
\newcommand{\inner}[2]{\langle #1,#2\rangle}
\title{An It\={o}--Wentzell Formula for Smooth Random Fields on Hilbert Spaces}
\author{Sascha Holl\\\texttt{sascha.holl@gmail.com}}
\date{Technical report --- November 2026}
\begin{document}
\maketitle
\begin{abstract}
We give a sufficient-condition version of the It\={o}--Wentzell formula for a random field whose spatial argument and range are separable Hilbert spaces and whose noise is cylindrical. The mixed spatial--noise correction is treated as a contraction of a Hilbert--Schmidt coefficient with a square-summable family of derivative operators. Spatial Galerkin projection reduces the proof to the classical finite-dimensional It\={o}--Wentzell formula with countably many Brownian motions. We also identify two obstructions to an unrestricted extension: smoothness of the coefficient into the Hilbert--Schmidt class alone does not control the mixed trace, and pointwise stochastic integration need not produce a continuous random field on an infinite-dimensional space. The purpose of this report is to provide a self-contained formulation and proof under explicit sufficient conditions.
\end{abstract}
\noindent\textbf{Keywords:} It\={o}--Wentzell formula; random fields; Hilbert spaces; cylindrical Wiener process; stochastic calculus.\par
\noindent\textbf{MSC 2020:} 60H05, 60H10, 60H15.

\section{Setting and motivation}
The It\={o}--Wentzell formula describes the composition of a semimartingale field with a semimartingale argument. Classical accounts include Kunita \cite{kunita}, and Krylov \cite[Theorem 3.1]{krylov} treats real-valued fields on Euclidean spaces driven by countably many Brownian motions under weaker spatial regularity. An especially relevant earlier work is Grecksch \cite{grecksch}, whose abstract describes an It\={o} formula for two-parameter Hilbert-space-valued random fields and results on martingale fields. The relationship between the hypotheses of that work and those used here is not addressed. We consider both an infinite-dimensional spatial argument and an infinite-dimensional range. The approach is deliberately restricted to fields for which a sufficiently regular version is already available. It does not assert existence of such a version from pointwise stochastic integrability.

Let $H,K,U$ be separable real Hilbert spaces and let $(\Omega,\mathcal A,\Pp)\text{ with filtration }(\mathcal F_t)_{t\in[0,T]}$ satisfy the usual conditions. Let $W$ be a cylindrical Wiener process on $U$ and $(e_j)$ an orthonormal basis of $U$; write $\beta_j=W(e_j)$. We use $\HS(E,F)$ for Hilbert--Schmidt operators and $\Lop(E,F)$ for bounded linear operators. For a bounded bilinear map $A:H\times H\to K$ and $B\in\HS(U,H)$ define
\begin{equation}\label{eq:quad}
 \mathcal T(A,B)=\sum_j A(Be_j,Be_j).
\end{equation}
This series is absolutely convergent, with norm at most $\norm{A}\norm{B}_2^2$, and independent of the orthonormal basis (by approximation of $B$ by finite-rank operators in Hilbert--Schmidt norm, or the projective tensor product).

\section{The mixed contraction}
\begin{proposition}\label{prop:contraction}
Suppose $D\in\Lop(H,\HS(U,K))$. For a fixed orthonormal basis $(e_j)$ of $U$, set $D_j v=(Dv)e_j$ and assume
\begin{equation}\label{eq:summable}
 \norm{D}_{\mathcal D,e}^2:=\sum_j\norm{D_j}_{\Lop(H,K)}^2<\infty.
\end{equation}
Then for every $B\in\HS(U,H)$ the sum
\begin{equation}\label{eq:mixed}
 \mathcal C(D,B)=\sum_j (D(Be_j))e_j
\end{equation}
is absolutely convergent and $\norm{\mathcal C(D,B)}\leq\norm{D}_{\mathcal D,e}\norm{B}_2$. Its value is independent of the basis among bases satisfying \eqref{eq:summable}; more generally its partial sums in any other basis converge to the same value, without asserting absolute convergence there.
A basis-independent sufficient condition for \eqref{eq:summable} in every basis is $D\in\HS(H,\HS(U,K))$.
\end{proposition}
\begin{proof}
The estimate follows from Cauchy--Schwarz. For a finite-rank operator $B=\sum_{i=1}^r v_i\otimes u_i$ (with $(v\otimes u)z=v\inner{u}{z}$), the sum equals $\sum_{i=1}^r(Dv_i)u_i$ in any basis. To justify convergence in any other basis, let $P$ be a finite-rank orthogonal projection on $U$. The corresponding contraction of $BP$ equals the sum over any orthonormal basis of $\operatorname{ran}P$ by the finite-rank identity. If $P_n\to I$ strongly, $BP_n\to B$ in Hilbert--Schmidt norm; the estimate in the original basis implies convergence to \eqref{eq:mixed}. Finally, for an orthonormal basis $(f_i)$ of $H$,
\[
 \sum_j\norm{D_j}^2\leq\sum_{j,i}\norm{(Df_i)e_j}^2=\norm{D}_{\HS(H,\HS(U,K))}^2.
\]
\end{proof}
\begin{remark}[Trace-class noise]
If $Q$ is a positive trace-class operator on a noise space $V$ and a coefficient $h:H\to\Lop(V,K)$ is Fr\'echet differentiable in operator norm, its effective coefficient $G(x)=h(x)Q^{1/2}$ satisfies \eqref{eq:summable}: $\norm{D_j}\leq\norm{Dh}\norm{Q^{1/2}e_j}$. Thus the hypothesis is compatible with the standard $Q$-Wiener representation. For an arbitrary coefficient defined only on the covariance space $V_0=Q^{1/2}V$ with its Cameron--Martin norm, this argument is not automatic. The effective coefficient is then $h(x)Q^{1/2}$, with $Q^{1/2}:V\to V_0$.
\end{remark}

\section{A Hilbert-space It\={o}--Wentzell formula}
\begin{theorem}\label{thm:main}
Let $x_0\in H$ be deterministic and let
\[
 X_t=x_0+\int_0^t b_s\dd s+\int_0^t B_s\dd W_s,
\]
where $b$ and $B$ are predictable, with values in $H$ and $\HS(U,H)$, respectively, and $\int_0^T(\norm{b_s}+\norm{B_s}_2^2)\dd s<\infty$ almost surely.
Let $F:\Omega\times[0,T]\times H\to K$ and predictable-in-time, jointly measurable fields $a$ and $G$ with values in $K$ and $\HS(U,K)$ satisfy, for each deterministic $x\in H$, indistinguishably in $t$,
\begin{equation}\label{eq:field}
 F_t(x)=F_0(x)+\int_0^t a_s(x)\dd s+\int_0^t G_s(x)\dd W_s.
\end{equation}
Assume the following.
\begin{enumerate}
\item $F(\cdot,x)$ is adapted for every $x$. On a single event of full probability, $F$ is jointly continuous in $(t,x)$, $F_t$ is twice continuously Fr\'echet differentiable, and its derivatives $DF_t,D^2F_t$ are jointly continuous in $(t,x)$ in operator norm. The fields and their derivatives are jointly measurable in the sense that all evaluations used below are predictable at predictable $H$-valued arguments.
\item For almost every $(\omega,t)$, $a_t$ is continuous and $G_t$ is continuously Fr\'echet differentiable into $\HS(U,K)$. The map $x\mapsto DG_t(x)$ satisfies \eqref{eq:summable} and is continuous for the norm $\norm{\cdot}_{\mathcal D,e}$ on that space. These fields and derivative evaluations are jointly predictable and Borel measurable.
\item For every $R<\infty$, the following bounds hold almost surely:
\begin{align*}
 &\sup_{t\leq T,\,\norm{x}\leq R}(\norm{F_t(x)}+\norm{DF_t(x)}+\norm{D^2F_t(x)})<\infty,\\
 &\int_0^T\sup_{\norm{x}\leq R}\norm{a_s(x)}\dd s<\infty,\\
 &\int_0^T\sup_{\norm{x}\leq R}\big(\norm{G_s(x)}_2^2+\norm{DG_s(x)}_{\mathcal D,e}^2\big)\dd s<\infty.
\end{align*}
\end{enumerate}
Then, almost surely for all $t\in[0,T]$,
\begin{align}\label{eq:main}
 F_t(X_t)-F_0(x_0)={}&\int_0^t\bigg(a_s(X_s)+DF_s(X_s)b_s\notag\\
 &\quad+\frac12\mathcal T(D^2F_s(X_s),B_s)+\mathcal C(DG_s(X_s),B_s)\bigg)\dd s\notag\\
 &+\int_0^t\big(G_s(X_s)+DF_s(X_s)B_s\big)\dd W_s.
\end{align}
In particular the two second-order corrections converge absolutely in $K$, and the equation is independent of the chosen orthonormal basis.
\end{theorem}
\begin{proof}
\emph{Well-definedness.} Continuous paths of $X$ are bounded on $[0,T]$. On the event that their bound is at most $R$, the assumptions, \eqref{eq:quad} and Proposition~\ref{prop:contraction} imply integrability of all drift terms. For the mixed correction use Cauchy--Schwarz in time. The diffusion term has finite time integral of its squared Hilbert--Schmidt norm. The stated measurability hypotheses and the countable sums ensure predictable integrands. All integrals are thus local stochastic or Bochner integrals.

\emph{Finite-dimensional spatial reduction.} Let $P_n$ project onto the span $H_n$ of the first $n$ vectors of an orthonormal basis of $H$, and put $X_t^n=x_0+P_n(X_t-x_0)$ and $Y_t^n=P_n(X_t-x_0)$. For a fixed $k\in K$, restrict $x\mapsto\inner{F_t(x_0+x)}{k}$ to $H_n\simeq\mathbb R^n$. This is a scalar semimartingale field with coefficients given by the corresponding restrictions of $a,G$. The spatial argument $Y^n$ starts at zero and has coefficients $P_nb,P_nB$. Apply the scalar It\={o}--Wentzell theorem with countably many driving Brownian motions \cite[Theorem 3.1]{krylov}. To see its hypotheses, bounded spatial balls in $H_n$ map into bounded balls of $H$; the drift and diffusion bounds give the required local time integrals. Spatial continuity of the contracted derivatives follows from the continuity in operator norm and \eqref{eq:summable}, including continuity of the sum of absolute mixed terms by Cauchy--Schwarz. The pathwise boundedness of $Y^n$ and of the field on bounded balls gives the localized space--time integrability conditions of that theorem. Equation \eqref{eq:field} supplies its pointwise semimartingale representation. Hence \eqref{eq:main} holds for the scalar pairing, with $X$ replaced by $X^n$, $b$ by $P_nb$ and $B$ by $P_nB$. All terms on the right define $K$-valued continuous processes; taking a countable dense set of $k$ gives the $K$-valued equality simultaneously in $t$.

\emph{Passage to the limit.} The image of each continuous path $X-x_0$ is compact, and strong convergence of the uniformly bounded projections is uniform on compact subsets. Thus $X^n\to X$ uniformly in time almost surely. Moreover $P_nb_s\to b_s$ and $P_nB_s\to B_s$ in the relevant norms for almost every $(\omega,s)$, the latter by Hilbert--Schmidt summability. All $X^n$ and $X$ lie within one path-dependent bounded ball. The continuity of $F$ on the compact closure of the union of the graphs of $X^n$ and $X$ yields uniform convergence of the left side. For the first two drift terms, continuity and the stated bounds give convergence in $L^1(0,T;K)$ almost surely. For the quadratic term use the bound
\begin{align*}
 \norm{\mathcal T(A_n,B_n)-\mathcal T(A,B)}
 &\leq \norm{A_n-A}\norm{B_n}_2^2\\
 &\quad+\norm{A}(\norm{B_n}_2+\norm{B}_2)\norm{B_n-B}_2.
\end{align*}
For the mixed term use the analogous bilinear estimate in the $\mathcal D,e$ and Hilbert--Schmidt norms. This gives convergence of both in $L^1$ by dominated convergence. The diffusion coefficients converge in $L^2(0,T;\HS(U,K))$ almost surely by the same argument. Convergence of the corresponding stochastic integrals uniformly in probability follows from the martingale maximal inequality after stopping at a fixed bound on the integrated squared difference. The deterministic integrals converge uniformly almost surely. Consequently the equality passes to the limit in the space of continuous processes with the uniform-in-probability topology, proving indistinguishability. Basis independence follows from the contractions and uniqueness of the stochastic integral.
\end{proof}

\section{Obstructions to weaker hypotheses}
\begin{proposition}[Failure of an unrestricted mixed trace]\label{prop:tracefail}
There exist $D\in\Lop(H,\HS(U,\mathbb R))$ and $B\in\HS(U,H)$ such that the formal series \eqref{eq:mixed} diverges, even when $D$ is the derivative at the origin of a globally bounded smooth $\HS(U,\mathbb R)$-valued function with bounded derivatives of every order.
\end{proposition}
\begin{proof}
Take $H=U=\ell^2$, $Be_j=j^{-1}e_j$, and let $Dv$ be the functional $u\mapsto\inner{v}{u}$. Then $B$ is Hilbert--Schmidt, and $D$ is an isometry into $\HS(U,\mathbb R)$. Nevertheless $(D(Be_j))e_j=j^{-1}$ and the harmonic series diverges. Choose a scalar $C^\infty$ cutoff $\chi$ with compact support and equal to one near zero, and put $G(x)u=\chi(\norm{x}^2)\inner{x}{u}$. Then $DG(0)=D$ and all derivatives of $G$ are bounded. This example concerns the coefficient contraction, not a random field satisfying Theorem~\ref{thm:main}.
\end{proof}
\begin{proposition}[Failure of a continuous modification]\label{prop:version}
Let $H=U=\ell^2$ and let $G$ be the smooth bounded coefficient just constructed. For every fixed $x$ the scalar field $F_t(x)=\int_0^t G(x)\dd W_s$ is well-defined, but for every fixed $t>0$ it admits no modification continuous at $0$ almost surely.
\end{proposition}
\begin{proof}
For $n\geq 2$ let $x_n=e_n/\sqrt{\log n}$. Then $x_n\to 0$ and eventually $F_t(x_n)=\beta_n(t)/\sqrt{\log n}$. For any $c\in(0,\sqrt{2t})$, independence and the Gaussian tail estimate give
\[
 \sum_n\Pp\big(|\beta_n(t)|>c\sqrt{\log n}\big)=\infty.
\]
The second Borel--Cantelli lemma implies that $|F_t(x_n)|>c$ infinitely often almost surely. Any modification must agree at the countably many deterministic points $0,x_2,x_3,\ldots$ on one event of full probability, so continuity at zero is impossible. The field can be realized in the original covariance-space notation for any injective trace-class $Q$ on $\ell^2$: if $U_0=Q^{1/2}U$, define $h(x):U_0\to\mathbb R$ by $h(x)v=G(x)Q^{-1/2}v$. Then $h(x)$ is Hilbert--Schmidt with respect to the covariance-space norm, and its $Q$-Wiener stochastic integral agrees with the cylindrical one.
\end{proof}
\begin{example}[The mixed term has coefficient one]
For real $X_t=\beta_t$ and $F_t(x)=x\beta_t$, the formal formula gives $\dd(\beta_t^2)=2\beta_t\dd\beta_t+\dd t$. Thus the mixed correction is not multiplied by $1/2$. On every bounded ball the field and its derivatives satisfy the local bounds in the theorem.
\end{example}

\section{Scope}
The formula assumes a regular version of the random field, boundedness on bounded sets, and quantitative control of its spatial noise derivative. It does not derive a regular version from local characteristics or address weak or mild solutions of stochastic evolution equations. The Galerkin argument uses the finite-dimensional theorem of Krylov \cite{krylov}; Hilbert-space random-field calculus has also been studied by Grecksch \cite{grecksch}. No assertion of priority or optimality of the hypotheses is made.


\begin{thebibliography}{9}
\bibitem{kunita} H. Kunita, \emph{Stochastic Flows and Stochastic Differential Equations}, Cambridge University Press, 1990.
\bibitem{krylov} N. V. Krylov, On the It\={o}--Wentzell formula for distribution-valued processes and related topics, \emph{Probability Theory and Related Fields} \textbf{150} (2011), 295--319. \href{https://doi.org/10.1007/s00440-010-0275-x}{doi:10.1007/s00440-010-0275-x}; \href{https://arxiv.org/abs/0904.2752}{arXiv:0904.2752}.
\bibitem{grecksch} W. Grecksch, Eine ITO-Formel f\"ur HILBERTraum-wertige zuf\"allige Felder, \emph{Mathematische Nachrichten} \textbf{101} (1981), 179--197. \href{https://doi.org/10.1002/mana.19811010116}{doi:10.1002/mana.19811010116}.
\bibitem{dpz} G. Da Prato and J. Zabczyk, \emph{Stochastic Equations in Infinite Dimensions}, 2nd ed., Cambridge University Press, 2014.
\end{thebibliography}
\end{document}